\documentclass[11pt]{amsart}
\usepackage[utf8]{inputenc}
\usepackage{amsthm, amsfonts, amssymb}
\usepackage{mathrsfs}
\usepackage{graphicx}
  \newcommand{\C}{\ensuremath{\mathbb{C}}}%
  \newcommand{\R}{\ensuremath{\mathbb{R}}}%
  \newcommand{\N}{\ensuremath{\mathbb{N}}}%
  \newcommand{\Z}{\ensuremath{\mathbb{Z}}}%
  \newcommand{\egdef}{\ensuremath{\stackrel{\text{\tiny def}}=}}%

\newtheorem{thm}{Theorem}[section]
\newtheorem{lemma}[thm]{Lemma}
\theoremstyle{remark}
\newtheorem*{rem}{Remark}

\author{Mikael de la Salle}
\title{A shorter proof of a result by Potapov and Sukochev on Lipschtiz functions on $S^p$}
\address{D\'epartement de Math\'ematiques et Applications \\ \'Ecole Normale
  Sup\'erieure  \\ 45 rue d'Ulm \\ 75005 Paris}

\address{Institut de Math\'ematiques de Jussieu \\ rue du Chevalleret \\
  75013 Paris}
\thanks{Partially supported by  ANR-06-BLAN-0015}
\email{mikael.de.la.salle@ens.fr}

\begin{document}
\maketitle
\begin{abstract}	
 In this short note we give a short proof of a recent result by Potapov and Sukochev (arXiv:0904.4095v1), stating that a Lipschitz function on the real line remains Lipschitz on the (self-adjoint part of) non-commutative $L_p$ spaces with $1<p<\infty$.
\end{abstract}

In the preprint \cite{PSArx} Potapov and Sukochev presented a proof of the following result, answering an open question going back at least to Krein (see also \cite{MR725454}):
\begin{thm}[Potapov and Sukochev]
\label{thm_potapov_sukochev}
 For any $1<p<\infty$ there exists a constant $C_p$ such that for any $1$-Lipschitz function $f:\R\to \R$
 and any two selfadjoint operators $A,B$ on $\ell^2$ with $A-B$ belonging to $p$-Schatten class $S^p$
 (\emph{i.e.} such that $\|A-B\|_p \egdef Tr(|A-B|^p)^{1/p}<\infty$), $f(A)-f(B) \in S^p$ and moreover
\[ \|f(A)-f(B)\|_p \leq C_p \|A - B\|_p.\]
\end{thm}

In this note we present a shortcut to their proof. Our proof still uses the two main ingredients of their proof, Lemma \ref{thm=multiplicateur_borne} and Lemma \ref{thm=Fourier} below, but it does not require an extrapolation type argument or the use of a factorization theorem for non-commutative $H^p$ spaces. Here we will restrict ourselves to the case when the operators $A$ and $B$ belong to a finite-dimensional matrix algebra, but the proof easily extends to the case when $A$ and $B$ belong to any von Neumann algebra with a normal trace, replacing the norm in $S^p$ by the norm in the corresponding non-commutative $L^p$-space. We will denote $S^p_n$ the space $M_n(\C)$ equipped the norm $\|A\|_p =Tr(|A-B|^p)^{1/p}$.

As explained in \cite{PSArx} (and also in Widom's problem 4.21 ``When are differentiable functions differentiable?'' of \cite{MR734178} or problem 6.4 of \cite{MR1334345}), the fact that a function $f$ is Lipschitz on $S^p_n$ is equivalent to the boundedness on $S^p_n$ of a class a Schur multiplier, namely the class of Schur multipliers with symbol $(f(\lambda_k)-f(\lambda_l))/(\lambda_k-\lambda_l) 1_{k \neq l}$ for any increasing sequence $\lambda_1<\dots<\lambda_n$ (more precisely the Lipschitz constant of $f$ on $S^p_n$ is equal to the supremum of the norms of these multipliers). Recall that for a family $\phi_{k,l} \in \C$ for $1\leq k,l\leq n$ we denote by $M_\phi$ and call ``Schur multiplier with symbol $\phi$'' the linear operator $M_n(\C) \to M_n(\C)$ sending a matrix $a=(a_{k,l})$ to the matrix $M_\varphi(a) = (\phi_{k,l} a_{k,l})$.

Since any $1$-Lipschitz function $f:\R \to \R$ is the difference of two \emph{non-decreasing} $1$-Lipschitz functions, it is also enough to treat the case when $f$ is non-decreasing (and even strictly increasing by an aproximation argument). Theorem \ref{thm_potapov_sukochev} is thus equivalent to the following theorem:
\begin{thm}
\label{thm=principal}
Let $1<p<\infty$. There is a constant $C_p>0$ such that for any (strictly) increasing $1$-Lipschitz function $f:\R \to \R$, any integer $n$ and any increasing sequence or real numbers $\lambda_1<\dots<\lambda_n$, 
if $\phi_{k,l} =(f(\lambda_k)-f(\lambda_l))/(\lambda_k-\lambda_l) 1_{k \neq l}$ then
\[ \left\| M_\phi \right\|_{S^p_n\to S^p_n} \leq C_p.\]
\end{thm}
\begin{rem}
 In fact both Potapov and Sukochev's proof and the following proof give the same bound for the completely bounded norm of $M_\phi$. This does not immediately follow from the statement of the Theorem since the question whether the boundedness of a Schur multiplier on $S^p$ implies its complete boundedness is still open when $1<p \neq 2<\infty$.
\end{rem}

Here are the two main ingredients from \cite{PSArx} we will use. The first one is a Fourier-transform trick (note that this kind of trick was already used in Lemma 1.7 of \cite{MR2372148}):
\begin{lemma}[\cite{PSArx}, Lemma 6]
\label{thm=Fourier}
There exists a function $g:\R\to\C$ such that:
\begin{itemize}
 \item $\int_\R |s|^m|g(s)| ds <\infty$ for any $m \in \N$.
\item	for any $0<\lambda<\mu$ we have
\[\frac{\lambda}{\mu} = \int_\R g(s) \lambda^{is} \mu^{-is} ds.\]
\end{itemize}
\end{lemma}
The second ingredient is the following, which is a consequence of the vector-valued Marcinkievicz multiplier theorem, due to Bourgain:
\begin{lemma}[\cite{PSArx}, Lemma 5]
\label{thm=multiplicateur_borne}
Let $1<p<\infty$. There exists $K_p >0$ such that for $s \in \R$, $n \in \N$ and if $M(s)$ the Schur multiplier on $M_n$ with symbol
$|k-l|^{is}$ (with the convention $0^{is}=0$). Then
\[\|M(s)\|_{S^p_n\to S^p_n} \leq  K_p (1+|s|).\]
\end{lemma}
\begin{rem}
In \cite{PSArx} this lemma is stated for the Schur multiplier $|k-l|^{is} 1_{k>l}$, but since $|k-l|^{is} = |k-l|^{is}1_{k<l}+ |k-l|^{is}1_{k>l}$ the above version follows from it.
\end{rem}

Let us now prove the main theorem:
\begin{proof}[Proof of Theorem \ref{thm=principal}]
Denote $\phi_{k,l} =(f(\lambda_k)-f(\lambda_l))/(\lambda_k-\lambda_l) 1_{k \neq l}$. Note that for $k \neq l$ since $f$ is increasing and $1$-Lipschitz we have $0 \leq \phi_{k,l} \leq 1$, and hence Lemma \ref{thm=Fourier} implies that
\[ \phi_{k,l} = \int_\R g(s) |f(\lambda_k)-f(\lambda_l)|^{is} |\lambda_k-\lambda_l|^{-is} ds.\]
For any increasing sequence $\mu_1<\dots<\mu_n$ denote by $M(s,(\mu_k)_k)$ the Schur multiplier with symbol $|\mu_k-\mu_l|^{is}$, so that
\[M_\phi = \int_\R g(s) M(s,(f(\lambda_k))_k) M(-s,(\lambda_k)_k) ds.\]
Since $\int_\R |g(s)| (1+|s|)^2 ds <\infty$ if we prove that $M(s,(\mu_k)_k)$ is bounded on $S^p_n$ with norm less than $K_p (1+|s|)$ (with $K_p$ given by Lemma \ref{thm=multiplicateur_borne}) then we will be done since it would imply that
\[ \|M_\phi\|_{S^p_n \to S^p_n} \leq K_p^2 \int_\R |g(s)| (1+|s|)^2 ds.\]
By a density argument it is enough to prove the bound on the norm of $M(s,(\mu_k)_k)$ when $\mu_k$ are all rational numbers. But then if $N$ is an integer such that $N\mu_k \in \Z$, the equality $|\mu_k-\mu_l|^{is} = N^{-is} |N\mu_k-N\mu_l|^{is}$ implies that we can assume that $\mu_k \in \Z$ for all $k$ (end even that $\mu_k \in \N$ by adding to the $\mu_k$'s a large enough number).

If $1\leq \mu_1<\dots<\mu_n$ then the matrix $(|\mu_k -\mu_l|^{is})_{1\leq k,l\leq n}$ is a submatrix of the matrix $(|k-l|^{is})_{1\leq k,l\leq \mu_n}$. The Schur multiplier $M(s,(\mu_k)_{1\leq k \leq n})$ is thus a restriction of the Schur multiplier $M(s,(k)_{1\leq k \leq \mu_n})$, which is just the multiplier $M(s)$ of Lemma \ref{thm=multiplicateur_borne} with $n$ replaced by $\mu_n$ and which is therefore bounded by $K_p (1+|s|)$.
\end{proof}
\bibliographystyle{plain}
\bibliography{biblio}

\end{document}